\documentclass[12pt]{amsart}
\newtheorem{theorem}{Theorem}
\newtheorem{proposition}[theorem]{Proposition}

\def\ds{\displaystyle}
\def\Re{\operatorname{Re}}

\title[An example of a bounded $\Bbb C$-convex domain]
{An example of a bounded $\Bbb C$-convex domain which is not
biholomorphic to a convex domain}

\author{Nikolai Nikolov, Peter Pflug and W\l odzimierz Zwonek}

\address
{Institute of Mathematics and Informatics\\ Bulgarian Academy of
Sciences\\ 1113 Sofia, Bulgaria}\email{nik@math.bas.bg}

\address{Carl von Ossietzky Universit\"at Oldenburg\\
Fachbereich Mathe\-ma\-tik\\ Postfach 2503\\ D-26111 Oldenburg,
Germany}\email{pflug@mathematik.uni-oldenburg.de}

\address{Instytut Matematyki, Uniwersytet Jagiello\'nski, Reymonta 4,
30-059 Krak\'ow, Poland}\email{Wlodzimierz.Zwonek@im.uj.edu.pl}
\begin{document}



\begin{thanks}{The paper was initiated while the third author's research
stay at the Carl von Ossietzky Universit\"at Oldenburg which was
supported by  the Alexander von Humboldt Foundation. The
third-named author was also supported by the KBN Research Grant
No. 1 PO3A 005 28. The first and second named authors were
supported by grants from DFG (DFG-Projekt 227/8).}
\end{thanks}

\begin{abstract} We show that the symmetrized bidisc is a $\Bbb
C$-convex domain. This provides an example of a bounded $\Bbb
C$-convex domain which cannot be exhausted by domains
biholomorphic to convex domains.
\end{abstract}

\maketitle

\section{Introduction}

Recall that a domain $D$ in $\Bbb C^n$ is called $\Bbb C$-convex
if any non-empty intersection with a complex line is contractible
(cf. \cite{APS, Hor}). A consequence of the fundamental Lempert
theorem (see \cite{Lem}) is the fact that any bounded $\Bbb C$-convex domain $D$ with
$C^2$ boundary has the following property (see \cite{Jac}):

{\it $(\ast)$ The Carath\'eodory distance and Lempert function of $D$
coincide.}

Any convex domain can be exhausted by smooth bounded convex ones (which are
obviously $\Bbb C$-convex); therefore, any convex domain satisfies $(\ast)$,
too. To extend this phenomenon to bounded $\Bbb C$-convex domains (see
Problem 4' in \cite{Zna}), it is sufficient to give a positive answer to one
of the following questions:

(a) {\it Can any bounded $\Bbb C$-convex domain be exhausted by
$C^2$-smooth $\Bbb C$-convex domains?} (See Problem 2 in
\cite{Zna} and Remark 2.5.20 in \cite {APS}.)

(b) {\it Is any bounded $\Bbb C$-convex domain biholomorphic to a
convex domain?} (See Problem 4 in \cite{Zna}.)

The main aim of this note is to give a negative answer to Question
(b).

Denote by $\Bbb G_2$ the so-called symmetrized bidisc, that is, the image of
the bidisc under the mapping whose components are the two elementary
symmetric functions of two complex variables. $\Bbb G_2$ serves as the first
example of a bounded pseudoconvex domain in $\Bbb C^2$ with the property
$(\ast)$ which cannot be exhausted by domains biholomorphic to convex
domains (see \cite{Cos,Edi}). We shall show that $\Bbb G_2$ is a $\Bbb
C$-convex domain. This fact  gives a counterexample to the question $(b)$
and simultaneously, it supports the conjecture that (cf. Problem 4' in
\cite{Zna}) {\it any bounded $\Bbb C$-convex domain has property $(\ast).$}
Note that the answer to the problem $(a)$ for $\Bbb G_2$ is not known. The
positive answer to this question would imply an alternative (to that of
\cite{Cos2} and \cite{AY}) proof of the equality of the Carath\'eodory
distance and Lempert function on $\Bbb G_2$ whereas the negative answer
would solve Problem 2 in \cite{Zna}.

Some additional properties of $\Bbb C$-convex domains and
symmetrized polydiscs are also given in the paper.

\section{Background and results}

Recall that a domain $D$ in $\Bbb C^n$ is called (cf.
\cite{Hor,APS}):

\begin{itemize}
\item {\it $\Bbb C$-convex} if any non-empty intersection with a
complex line is contractible (i.e. $D\cap L$ is connected and simply connected for any
complex affine line $L$
such that $L\cap D$ is not empty);

\item {\it linearly convex} if its complement in $\Bbb C^n$ is a
union of affine complex hyperplanes;

\item {\it weakly linearly convex} if for any $a\in\partial D$
there exists an affine complex hyperplane through $a$ which does
not intersect D.
\end{itemize}

Note that the following implications hold

\centerline{$\Bbb C$-convexity $\Rightarrow$ linear convexity
$\Rightarrow$ weak linear convexity.}

Moreover, these three notions coincide in the case of bounded
domains with $C^1$ boundary.

Let $\Bbb D$ denote the unit disc in $\Bbb C.$ Let
$\pi_n=(\pi_{n,1},\dots,\pi_{n,n}):\Bbb C^n\to\Bbb C^n$ be defined as
follows: $$ \pi_{n,k}(\mu)=\sum_{1\le j_1<\dots<j_k\le
n}\mu_{j_1}\dots,\mu_{j_k},\ \ 1\le k\le n, \;
\mu=(\mu_1,\dots,\mu_n)\in\Bbb C^n. $$ The set $\Bbb G_n:=\pi_n(\Bbb D^n)$
is called {\it the symmetrized n-disc } (cf. \cite{AY}, \cite{Jar-Pfl}).

Recall that $\Bbb G_2$ is the first example of a bounded pseudoconvex domain
with the property $(\ast)$ which cannot be exhausted by domains
biholomorphic to convex ones (see \cite{Cos,Edi}). On the other hand, $\Bbb
G_n,$ $n\ge 3,$ does not satisfy the property $(\ast)$ (see \cite{NPZ}). In
particular, it cannot be exhausted by domains biholomorphic to convex
domains, either.

In this note we shall show the following additional
properties of domains $\Bbb
G_n,$ $n\ge 2.$

\begin{theorem} {\rm (i)} $\Bbb G_2$ is a $\Bbb C$-convex domain.

{\rm (ii)} $\Bbb G_n,$ $n\ge 3,$ is a linearly convex domain which is
not $\Bbb C$-convex.
\end{theorem}

Theorem 1 (i) together with a result of \cite{Cos} and
\cite{Edi} gives a negative answer to the following question posed by S. V.
Znamenski\u{i} (cf. Problem 4 in \cite{Zna}):

{\it Is any bounded $\Bbb C$-convex domain biholomorphic to a
convex domain?}

Moreover, it seems to us that Theorem 1 (ii) gives the
first example of a linearly convex domain homeomorphic to $\Bbb
C^n,$ $n\ge 3,$ which is not $\Bbb C$-convex, is not a Cartesian product and does not satisfy
property $(\ast).$ To see that $\Bbb G_n$ is homeomorphic to $\Bbb
C^n,$ observe that $\rho_\lambda(z):=(\lambda z_1,\lambda^2
z_2,\dots,\lambda^n z_n)\in\Bbb G_n$ if $z\in\Bbb G_n$ and
$\lambda\in\Bbb C.$ Then setting $\ds h(z)=\max_{1\le
j\le n}\{|\mu_j|:\pi_n(\mu)=z\}$ and $g(z)=\frac{1}{1-h(z)},$ it
is easy to see that the function $\Bbb G_n\owns z\mapsto\rho_{g(z)}(z)\in\Bbb C^n$ is the desired
homeomorphism.

These remarks also show that $\Bbb G_n$ is close, in some sense, to
a balanced domain, that is, a domain $D$ in $\Bbb C^n$ such that
$\lambda z\in D$ for any $z\in D$ and $\lambda\in\overline{\Bbb
D}.$ On the other hand, in spite of the properties of $\Bbb G_n,$
one has the following.

\begin{proposition} Any weakly linearly convex balanced domain is
convex.
\end{proposition}

This proposition is a simple extension of Example 2.2.4 in
\cite{APS}, where it is shown that any $\Bbb C$-convex complete
Reinhardt domain is convex.

We may also prove some general property of $\Bbb C$-convex domains
showing that all {\it non-degenerate} $\Bbb C$-convex domains,
that is, containing {\it no} complex lines, are $c$-finitely
compact. For definitions of the Carath\'eodory distance $c_D$ of
the domain $D$, $c$-finite compactness, $c$-completeness and basic
properties of these notions we refer the Reader to consult
\cite{JP1}.

Observe that a degenerate linearly convex domain $D$ is linearly
equivalent to $\Bbb C\times D'$ (cf. Proposition 4.6.11 in
\cite{Hor}). Indeed, we may assume that $D$ contains the
$z_1$-line. Since the complement $^cD$ of $D$ is a union of
complex hyperplanes disjoint from this line, then $^cD=\Bbb
C\times G$ and hence $D=\Bbb C\times{^cG}.$ On the other hand, we
have

\begin{proposition} Any non-degenerate $\Bbb C$-convex
domain is biholomorhic to a bounded domain and $c$-finitely
compact. In particular, it is $c$-complete and hyperconvex.
\end{proposition}

\noindent{\bf Remarks.} (i) In virtue of Proposition 3, we claim
that one may conjecture more than Question (a) (see \cite{ZZ}),
namely, any $\Bbb C$-convex domain containing no complex
hyperplanes can be exhausted by bounded $C^2$-smooth $\Bbb
C$-convex domains (this is not true in general without the above
assumption); then the Carath\'eodory pseudodistance and Lempert
function will coincide on any $\Bbb C$-convex domain. \vskip1mm

\noindent (ii) The hyperconvexity of $\Bbb G_n$ is simple and
well-known (see \cite{Edi-Zwo}). The above proposition implies
more in dimension two. Namely, it implies that the symmetrized
bidisc is $c$-finitely compact. Although the symmetrized polydiscs
in higher dimensions are not $\Bbb C$-convex the conclusion of the
above proposition, that is, the $c$-finite compactness of the
symmetrized $n$-disc $\Bbb G_n$, holds for any $n\geq 2$. In fact,
it is a straightforward consequence of Corollary 3.2 in
\cite{Cos3}. \vskip1mm

\noindent{(iii)} Finally, we mention that, for $n\geq 2$, $\Bbb
G_n$ is starlike with respect to the origin if and only if $n=2.$
This observation gives the next difference in the geometric shape
of the $2$-dimensional and higher dimensional symmetrized discs.
Recall that the fact that  $\Bbb G_2$ is starlike is contained in
\cite{AY}. For the converse just take the point
$(3,3,1,0,\dots,0)$.

\section{Proofs}

\noindent{\it Proof of Theorem 1 {\rm (i)}.} We shall make use of
the following description of $\Bbb C$-convex domains. For
$a\in\partial D$, denote by $\Gamma(a)$ the set of all hyperplanes
through $a$ and disjoint from $D$. Then a bounded domain $D$ in
$\Bbb C^n$, $n>1,$ is $\Bbb C$-convex if and only if any
$a\in\partial D$ the set $\Gamma(a)$ is non-empty and connected as
a set in $\Bbb C\Bbb P^n$ (cf. Theorem 2.5.2 in \cite{APS}).

So we have to check that $\Gamma(a)$ is non-empty and connected
for any $a\in\partial\Bbb G_2.$

Let us first consider a regular point $\partial\Bbb G_2$, that is,
a point of the form $\pi_2(\mu)$, where $|\mu_1|=1$, $|\mu_2|<1$
(or vice versa). Then the complex tangent line to $\partial D$ at
$a$ is of the form $\{\pi_2(\mu_1,\lambda) :\lambda\in\Bbb C\}$,
which is obviously disjoint from $\Bbb G_2$. So $\Gamma(a)$ is a
singleton.

Now we fix a non-regular point of $\partial\Bbb G_2$, that
is, a point of the form $\pi_2(\mu)$, where $|\mu_1|=|\mu_2|=1.$

After a rotation we may assume that $\mu_1\mu_2=1$, that is,
$\mu_2=\bar\mu_1$. Then $\mu_1+\mu_2=2\Re \mu_1=:2x$, where
$x\in[-1,1]$.

We shall find all the possible directions of complex lines passing
simultaneously through $\pi_2(\mu)$ and an element of $\Bbb G_2$.
Any such line is of the form $\pi_2(\mu)+\Bbb
C(\pi_2(\mu)-\pi_2(\lambda))$, where $\lambda\in\Bbb D^2$. So the
$$A:={^c\Gamma}(\pi_2(\mu))=\{\frac{\lambda_1+\lambda_2-2x}{\lambda_1\lambda_2-1}:
\lambda_1,\lambda_2\in\Bbb D\}.$$ In particular,
$\Gamma(\pi_2(\mu))\neq\emptyset.$

To show the connectedness of $\Gamma(\pi_2(\mu)),$ we shall check
the simple-connectedness of $A.$ Let us recall that the mapping
$\frac{z-\alpha}{z-\beta}$, where $|\beta|>1$, maps the unit disc
$\Bbb D$ into the disc
$\triangle(\frac{1-\alpha\bar\beta}{1-|\beta|^2},
\frac{|\alpha-\beta|}{|\beta|^2-1})$, so
$$\{\frac{\lambda+\lambda_1-2x}{\lambda\lambda_1-1}:\lambda\in\Bbb
D\}=\triangle(\frac{2x-2\Re\lambda_1}{1-|\lambda_1|^2},
\frac{|2x\lambda_1-\lambda_1^2-1|}{1-|\lambda_1|^2})=:A_{\lambda_1}.
$$ Consequently the set $A=\bigcup\sb{\lambda_1\in\Bbb
D}A_{\lambda_1}\subset\Bbb C$ is simply connected. \qed

\

\noindent{\it Proof of Theorem 1 {\rm (ii)}.} For the proof of the
linear convexity of $\Bbb G_n$ consider the point
$z=\pi_n(\lambda)\in\Bbb C^n\setminus\Bbb G_n$. We may assume that
$|\lambda_1|\geq 1$. Then the set $$
B:=\{\pi_n(\lambda_1,\mu_1,\ldots,\mu_{n-1}):\mu_1,\ldots,\mu_{n-1}\in\Bbb
C\} $$ is disjoint from $\Bbb G_n$. On the other hand, it is easy
to see that $$ B=\{(\lambda_1+z_1,\lambda_1z_1+z_2,
\ldots,\lambda_1z_{n-2}+z_{n-1},\lambda_1z_{n-1}):z_1,\ldots,z_{n-1}\in\Bbb
C\}, $$ so $B$ is a complex affine hyperplane. Hence $\Bbb G_n$ is
linearly convex.

To show that $\Bbb G_n$ is not $\Bbb C$-convex for $n\geq 3,$
consider the points
$$a_t:=\pi_n(t,t,t,0,\ldots,0)=(3t,3t^2,t^3,0,\ldots,0), $$
$$b_t:=\pi_n(-t,-t,-t,0,\ldots,0)=(-3t,3t^2,-t^3,0,\ldots,0),\;
t\in(0,1). $$ Obviously $a_t,b_t\in\Bbb G_n$. Denote by $L_t$ the
complex line passing through $a_t$ and $b_t,$ that is,
$$L_t=\{c_{t,\lambda}:=(3t(1-2\lambda),3t^2,t^3(1-2\lambda),0,\ldots,0):
\lambda\in\Bbb C\}.$$ Assume that the set $\Bbb G_n\cap L_t$ is
connected. Since $a_t=c_{t,0}$ and $b_t=c_{t,1},$ then
$c_{t,\lambda}\in\Bbb G_n$ for some $\lambda=\frac{1}{2}+i\tau,$
$\tau\in\Bbb R.$ It follows that $$c_{t,\lambda}=(-6i\tau
t,3t^2,-2i\tau t^3,0,\ldots,0).$$ We may choose $\mu\in\Bbb D^n$
such that $\mu_j=0$, $j=4,\ldots,n,$ and
$c_{t,\lambda}=\pi_n(\mu),$ $\mu\in\Bbb D^n.$ Then
$-36\tau^2t^2=(\mu_1+\mu_2+\mu_3)^2=\mu_1^2+\mu_2^2+\mu_3^2+6t^2$
and hence
$$t^2=\frac{|\mu_1^2+\mu_2^2+\mu_3^2|}{36\tau^2+6}<\frac{3}{36\tau^2+6}
\leq\frac{1}{2}.$$ Therefore, $\Bbb G_n\cap L_t$ is not connected
if $t\in[\frac{1}{\sqrt{2}},1)$ and so $\Bbb G_n$ is not a $\Bbb
C$-convex domain.\qed

\

\noindent{\it Proof of Proposition 2.} Set $D^\ast:=\{w\in\Bbb
C^n:<z,w>\neq 1,\ \forall z\in D\}.$ We shall use the fact that a
domain $D$ in $\Bbb C^n$ containing the origin is weakly linearly
convex if and only if $D$ is a connected component of
$D^{\ast\ast}$ (cf. Proposition 2.1.4 in \cite{APS}).

Since our domain $D$ is balanced, it is easy to see that $D^\ast$
is balanced. We shall show $D^\ast$ is convex. Then, applying this
fact to $D^\ast,$ we conclude that $D^{\ast\ast}$ is a convex
balanced domain. On the other hand, it follows by our assumption
that $D$ is a component of $D^{\ast\ast}$ and hence
$D^{\ast\ast}=D.$

To see that $D^\ast$ is convex, suppose the contrary. Then we find
points $w_1,w_2\in D^\ast,$ $z\in D$ and a number $t\in(0,1)$ such
that $<z,tw_1+(1-t)w_2>=1$. We may assume that $|<z,w_1>|\ge 1.$
Since $D$ is balanced, we get $\tilde z:=\frac{z}{<z,w_1>}\in D$
and $<\tilde z,w_1>=1,$ a contradiction.\qed

\

\noindent{\it Proof of Proposition 3.} Let $D$ be non-degenerate
$\Bbb C$-convex domain in $\Bbb C^n.$ For any point $z\in{^cD}$
consider a hyperplane $L_z$ through $z$ and disjoint from $D.$ Let
$l_z$ be the orthogonal line through $0$ and orthogonal to $L_z.$
Denote by $\pi_z$ the orthogonal projection of $\Bbb C^n$ onto
$l_z$ and set $a_z=\pi_z(a).$ Observe that $D_z=\pi_z(D)$ is
biholomorphic to $\Bbb D,$ since it is connected, simply connected
(cf. Theorem 2.3.6  in \cite{APS}) and $\pi_z(z)\not\in\pi_z(D).$
Moreover, since $D$ is a non-degenerate linearly convex domain, it
is easy to see that there are $n$ $\Bbb C$-independent $l_z'$s. We
may assume that these $l_z$ are the set $C$ of coordinate planes.
Then $\ds D\subset G:=\prod_{l_z\in C}\pi_z(D)$ and $G$ is
biholomorphic to the polydisc $\Bbb D^n.$ In particular, $D$ is
biholomorphic to a bounded domain, hence it is $c$-hyperbolic.

Further, we may assume that $0\in D.$ To see that $D$ is
$c$-finitely compact, it is enough to show that $\ds\lim_{a\to
z}c_D(0;a)=\infty$ for any $z\in\partial D$ and, if $D$ is
unbounded, $z=\infty.$ But the last one follows by the fact that
$G$ is $c$-finitely compact. On the other hand, if $a\to
z\in\partial D,$ then $a_z\to\pi_z(z)\in\partial D_z$ and hence
$c_D(0;a)\ge c_{D_z}(0;a_z)\to\infty.$\qed


\begin{thebibliography}{}

\bibitem{AY} J.~Agler, N.~J.~Young, {\it
The hyperbolic geometry of the symmetrized bidisc}, J. Geom. Anal.
14 (2004), 375--403.

\bibitem{APS} M.~Andersson, M.~Passare, R.~Sigurdsson, {\it Complex convexity and
analytic functionals,} Birkh\"auser, Basel--Boston--Berlin, 2004.

\bibitem{Cos} C.~Costara, {\it The symmetrized bidisc and Lempert's theorem},
Bull. London Math. Soc. 36 (2004), 656--662.

\bibitem{Cos2} C.~Costara, {\it Dissertation}, Universit\'e
Laval (2004).

\bibitem{Cos3} C.~Costara, {\it On the spectral Nevanlinna--Pick problem},
Studia Math. 170 (2005), 23--55.

\bibitem{Edi} A. Edigarian, {\it A note on Costara's paper},
Ann. Polon. Math. 83 (2004), 189--191.

\bibitem{Edi-Zwo} A.~Edigarian, W.~Zwonek,
{\it Geometry of the symmetrized polydisc}, Arch. Math. (Basel) 84
(2005), 364--374.

\bibitem{Jac} D. Jacquet, {\it $\Bbb C$-convex domains with $C^2$ boundary},
Complex Variables and Elliptic Equations 51 (2006), 303--312.

\bibitem{Hor} L.~H\"ormander, {\it Notions of convexity},
Birkh\"auser, Basel--Boston--Berlin, 1994.

\bibitem{JP1} M.~Jarnicki, P.~Pflug, {\it Invariant
Distances and Metrics in Complex Analysis}, Walter de Gruyter, 8 (1993).

\bibitem{Jar-Pfl} M.~Jarnicki, P.~Pflug, {\it Invariant
distances and metrics in complex analysis--revisited}, Diss. Math.
430 (2005), 1--192.

\bibitem{Lem} L.~Lempert, {\it La m\'etrique de Kobayashi et la
repr\'esentation des domaines sur la boule}, Bull. Soc. Math.
France 109 (1981), 427--474.

\bibitem{NPZ} N.~Nikolov, P.~Pflug, W.~Zwonek,
{\it The Lempert function of the symmetrized polydisc in higher
dimensions is not a distance}, Proc. Amer. Math. Soc., to appear
(arXiv:math.CV/0601367).

\bibitem{Zna} S.~V. Znamenski\u{i}, {\it Seven $\Bbb C$-convexity problems}
(in Russian), Complex analysis in modern mathematics. On the 80th
anniversary of the birth of Boris Vladimirovich Shabat, E.~M.
Chirka (ed.), FAZIS, Moscow, 2001, 123--131.

\bibitem{ZZ} S.~V. Znamenski\u{i}, L.~N.~Znamenskaya, {\it Spiral
connectedness of the sections and projections of $\Bbb C$-convex
sets}, Math. Notes 59 (1996), 253--260.

\end{thebibliography}
\end{document}